\documentclass{amsart}

\usepackage{xypic}
\input xy
\xyoption{all}
\usepackage{epsfig}
\usepackage{amsthm}
\usepackage{amssymb}
\usepackage{amsmath}
\usepackage{amscd}
\usepackage{xcolor}

\newcommand{\bg}{\begin{equation}}
\newcommand{\ed}{\end{equation}}
\newcommand{\bga}{\begin{eqnarray}}
\newcommand{\eda}{\end{eqnarray}}
\newcommand{\pf}{\textbf{Proof:\ }}

\def\cbdu{\par{\raggedleft$\Box$\par}}

\newtheorem {Theorem}  {Theorem}

\numberwithin{Theorem}{section}

\newtheorem {Lemma}[Theorem]  {Lemma}
\newtheorem {Proposition}[Theorem]{Proposition}
\theoremstyle{definition}
\newtheorem{Definition}[Theorem]{Definition}
\theoremstyle{remark}

\expandafter\chardef\csname pre amssym.def
at\endcsname=\the\catcode`\@ \catcode`\@=11
\def\undefine#1{\let#1\undefined}
\def\newsymbol#1#2#3#4#5{\let\next@\relax
 \ifnum#2=\@ne\let\next@\msafam@\else
 \ifnum#2=\tw@\let\next@\msbfam@\fi\fi
 \mathchardef#1="#3\next@#4#5}
\def\mathhexbox@#1#2#3{\relax
 \ifmmode\mathpalette{}{\m@th\mathchar"#1#2#3}%
 \else\leavevmode\hbox{$\m@th\mathchar"#1#2#3$}\fi}
\def\hexnumber@#1{\ifcase#1 0\or 1\or 2\or 3\or 4\or 5\or 6\or 7\or 8\or
 9\or A\or B\or C\or D\or E\or F\fi}

\font\teneufm=eufm10 \font\seveneufm=eufm7 \font\fiveeufm=eufm5
\newfam\eufmfam
\textfont\eufmfam=\teneufm \scriptfont\eufmfam=\seveneufm
\scriptscriptfont\eufmfam=\fiveeufm

\catcode`\@=\csname pre amssym.def at\endcsname

\newcounter{remark}
\renewcommand{\a}{\alpha}

\renewcommand{\b}{\beta}

\newcommand{\R}{\mathbf{R}}

\renewcommand{\div}{\mbox{div}}

\def  \R   {{\mathbb R}}

\def  \12  {{\frac{1}{2}}}

\def\build#1_#2^#3{\mathrel{\mathop{\kern 0pt#1}\limits_{#2}^{#3}}}

\begin{document}

\title[Well-posedness]{On well-posedness of generalized Hall-magneto-hydrodynamics}


\author [Mimi Dai]{Mimi Dai}
\address{Department of Mathematics, Stat. and Comp. Sci.,  University of Illinois Chicago, Chicago, IL 60607,USA}
\email{mdai@uic.edu} 

\author [Han Liu]{Han Liu}
\address{Department of Mathematics, Stat. and Comp. Sci.,  University of Illinois Chicago, Chicago, IL 60607,USA}
\email{hliu94@uic.edu} 

\thanks{The work of the authors was partially supported by NSF Grant DMS--1815069.}





\begin{abstract}

We obtain local well-posedness result for the generalized Hall-magneto-hydrodynamics system in Besov spaces ${\dot B^{-(2\alpha_1-\gamma)}_{\infty, \infty}} \times {\dot B^{-(2\alpha_2-\beta)}_{\infty, \infty}(\mathbb R^3)}$ with suitable indexes $\alpha_1, \alpha_2, \beta$ and $\gamma.$ As a corollary, the hyperdissipative electron magneto-hydrodynamics system is globally well-posed in ${\dot B^{-(2\alpha_2-2)}}_{\infty, \infty}(\mathbb R^3)$ for small initial data.

\bigskip

KEY WORDS: Well-posedness, Hall-MHD, Electron-MHD,  Small data.

\hspace{0.02cm}CLASSIFICATION CODE: 	35Q35, 35Q60, 35A05
\end{abstract}

\maketitle

\section{Introduction}

In this paper, we study the well-posedness problem of the following generalized Hall-magneto-hydrodynamics (Hall-MHD) system
\begin{equation}\label{fhmhd}
\begin{cases}
u_t+(u\cdot \nabla)u-(b \cdot \nabla)b+\nabla p = -\nu(-\Delta)^{\a_1} u,\\
b_t+(u \cdot \nabla )b-(b \cdot \nabla) u+ \eta \nabla \times ((\nabla \times b) \times b)= -\mu (-\Delta)^{\a_2} b,\\
\nabla \cdot u=0, \ \nabla \cdot b=0,\\
u(0,x)=u_0, \ b(0, x)=b_0, \ t \in \R^+,  x\in \R^3,
\end{cases}
\end{equation}
with the parameters $\a_1, \a_2 >0$ and the constants $\nu, \mu >0, \eta \geq 0.$ 

In particular, the fourth term on the left-hand side of the second equation is called the Hall term. When $\a_1=\a_2=1,$ $\eta>0,$ system (\ref{fhmhd}) becomes the standard Hall-MHD system, whereas the case $\eta=0$ corresponds to the generalized magneto-hydrodynamics (MHD) system. 

Derived in \cite{ADFL} as the incompressible limit of a two-fluid isothermal Euler-Maxwell system for electrons and ions, the Hall-MHD system describes the evolution of a system consisting of charged particles that can be approximated as a conducting fluid, in the presence of a magnetic field $b,$ with $u$ denoting the fluid velocity, $p$ the pressure, $\nu$ the viscosity, $\mu$ the magnetic resistivity and $\eta$ a constant  determined by the ion inertial length. The MHD and Hall-MHD systems have a wide range of applications in plasma physics and astrophysics, including modelling solar wind turbulence, designing tokamaks as well as studying the origin and dynamics of the terrestrial magnetosphere. Notably, the Hall-MHD system serves a vital role in interpreting the magnetic reconnection phenomenon, frequently observed in space plasmas. For more physical backgrounds, we refer readers to \cite{C, G1, G2, GB, L, PM}.

Over the past decade, various mathematical results concerning the Hall-MHD system have been obtained. A mathematically rigorous derivation of the system is due to Acheritogaray, Degond, Frouvelle and Liu \cite{ADFL}. Concerning the solvability of the system, Chae, Degond and Liu \cite{CDL} obtained global-in-time existence of weak solutions and local-in-time existence of classical solutions.  In \cite{CL}, Chae and Lee established a blow-up criterion and a small data global existence result. In addition, local well-posedness results can be found in the works by Dai \cite{D2, D3}, and global existence results for small data were also proved by Wan and Zhou \cite{WZ1} as well as by Kwak and Lkhagvasuren \cite{KL}. For various regularity criteria, readers are referred to \cite{D1, FFNZ, FLN, HAHZ, WL, Y3, Y4, YZ0, Z}.  Regarding the propeties of the solutions, the temporay decay of weak solutions was studied by Chae and Schonbek \cite{CS}, while the stability of global strong solutions is due to Benvenutti and Ferreira \cite{BF}. On the other hand, in the irresistive setting, there are striking ill-posedness results due to Chae and Weng \cite{CW} as well as Jeong and Oh \cite{JO}. Recently, Dai \cite{D4} proved the non-uniqueness of the Leray-Hopf weak solution via a convex integration scheme.   

The generalized system (\ref{fhmhd}) has also attracted mathematicians' attentions. Chae, Wan and Wu \cite{CWW} proved local well-posedness in the case $\a_1=0,$ $\a_2>\frac{1}{2},$ while local well-posedness result for $0<\a_1 \leq 2,$ $1 < \a_2 \leq 2$ and global well-posedness result for $\a_1 \geq \frac{5}{4},$  $\a_2 \geq \frac{7}{4}$ were obtained respectively by Wan and Zhou \cite{WZ2} and Wan \cite{W}. Small data global solutions were established in \cite{PMZ, WYT, Y1}. In addition, decay results of global smooth solutions in the cases where either $\a_1$ or $\a_2=0$ is due to Dai and Liu \cite{DL}. We refer readers to \cite{FSZ, GMS, JZ, PZ} for a number of regularity criteria.

In this paper, we shall prove that system (\ref{fhmhd}) is locally well-posed in the Besov space ${\dot B^{-(2\a_1-\gamma)}_{\infty, \infty}} \times {\dot B^{-(2\a_2-\b)}_{\infty, \infty}}(\R^3)$ for suitable choices of $\a_1, \a_2, \beta$ and $\gamma.$ Our main result states as follows.
\begin{Theorem}[Local well-posedness]\label{lwp}
For $(u_0, b_0) \in {\dot B^{-(2\a_1-\gamma)}_{\infty, \infty}} \times {\dot B^{-(2\a_2-\b)}_{\infty, \infty}}(\R^3),$ there exists a unique local-in-time solution $(u,b)$ to system (\ref{fhmhd}) such that $$(u,b) \in L^\infty\big(0,T; {\dot B^{-(2\a_1-\gamma)}_{\infty, \infty}} \times {\dot B^{-(2\a_2-\b)}_{\infty, \infty}(\R^3)}\big)$$ with $T = T \big(\nu, \mu, \eta, \|u_0\|_{\dot B^{-(2\a_1-\gamma)}_{\infty, \infty}}, \|b_0\|_{\dot B^{-(2\a_2-\b)}_{\infty, \infty}}\big),$ provided that the parameters $\a_1, \a_2, \b$ and $\gamma$ satisfy the following constraints 
\begin{equation}
\begin{cases}
\gamma \geq \max\{1, \frac{\a_1}{\a_2}\}, \\
\b \geq \max\{ 2, \frac{(\gamma +1)\a_2}{2\a_1}\},\\
\frac{\gamma}{2}< \a_1 < \gamma,\\
\frac{\b}{2} <\a_2<\b.
\end{cases} 
\end{equation}
\end{Theorem}

An interesting byproduct of the above result is small data global well-posedness for the electron MHD (EMHD) equations, the fluid-free version of system (\ref{fhmhd}).
\begin{Theorem}[Global existence for small data]\label{gwp}
Let $1<\a_2 <2.$ There exists some $\varepsilon=\varepsilon(\mu)>0$ such that if $\|b_0\|_{\dot B^{-(2\a_2-2)}_{\infty, \infty}(\R^3)} \leq \varepsilon,$ then there exists a solution $b$ to the EMHD system, i.e., system (\ref{fhmhd}) with $u \equiv 0,$ satisfying $$b \in L^\infty\big(0, +\infty; {\dot B^{-(2\a_2-2)}_{\infty, \infty}(\R^3)}\big) \text{ and } \sup_{t>0} t^{\frac{\a_2-1}{\a_2}}\|b\|_{L^\infty(\R^3)} < \infty.$$
\end{Theorem}

For generalized MHD system, local and global well-posedness results in Besov spaces were proved in \cite{YZ} via the same mechanism as the one in this paper, in spite of a major difference between the MHD and Hall-MHD systems in terms of scaling properties. In brief, the generalized MHD system scales as 
$$
u_\lambda(t,x)=\lambda^{2\a_1-1}u(\lambda^{2\a_1}t, \lambda x), \ b_\lambda(t,x)=\lambda^{2\a_2-1}b(\lambda^{2\a_2}t, \lambda x),
$$
while the EMHD equations scale as 
$ b_\lambda(t,x)=\lambda^{2\a_2-2}b(\lambda^{2\a_2}t, \lambda x),$ resulting in an absence of scaling invariance along with a lack of the notion of criticality in the Hall-MHD system, which seems to render the global well-posedness for the full system (\ref{fhmhd}) rather elusive. For system (\ref{fhmhd}), we can only establish local well-posedness, in contrast to the generalized MHD system, which possesses global-in-time solutions in the largest critical space $\dot B^{-(2\a_1-1)}_{\infty, \infty} \times \dot B^{-(2\a_2-1)}_{\infty, \infty}(\R^3),$ with $\a_1=\a_2, \ \frac{1}{2}< \a_1, \a_2 <1$ for small initial data, as proven in \cite{YZ}. The fact that the well-posedness result for the Hall-MHD system deviates from that for the MHD system is an evidence that the new scale and non-linear interactions introduced by the Hall term $\nabla \times ((\nabla \times b)\times b)$ play a significant role. 

\bigskip

\section{Preliminaries}
\label{sec:pre}

\subsection{Notation}
Throughout the paper, we will use $C$ to denote different constants. The notation $A \lesssim B$ means that $A \leq CB$ for some constant $C.$  For simplicity, we denote the caloric extensions $e^{-\nu t (-\Delta)^{\a_1}}u_0$ and $e^{-\mu t(-\Delta)^{\a_2}}b_0$ by $\tilde u_0$ and $\tilde b_0,$ respectively. In addition, we use $\mathbb{P}$ to denote the Helmholtz-Leray projection onto solenoidal vector fields, which acts on a vector field $\phi$ as $$\mathbb{P}\phi = \phi +\nabla \cdot (-\Delta)^{-1} \div \phi.$$ 

\subsection{Besov spaces via Littlewood-Paley theory}
\label{sec:LPD}
We shall briefly recall the homogeneous Littlewood-Paley decomposition, through which we shall define the homogeneous Besov space. For a complete description of Littlewood-Paley theory and its applications, we refer readers to \cite{BCD, G}.

We introduce the radial function $\chi\in C_0^\infty(\R^n)$ such that $0 \leq \chi \leq 1$ and
\begin{equation}\notag
\chi(\xi)=
\begin{cases}
1, \ \ \mbox { for } |\xi|\leq\frac{3}{4}\\
0, \ \ \mbox { for } |\xi|\geq 1.
\end{cases}
\end{equation}
Let $\varphi \in C_0^\infty(\R^n)$ be such that $\varphi(\xi)=\chi(\xi/2)-\chi(\xi).$ We construct a family of smooth functions $\{\varphi_q \}_{q \in \mathbb{Z}}$ supported on dyadic annuli in the frequency space, defined as 
\begin{equation}\notag
\varphi_q(\xi)=\varphi(2^{-q}\xi), \ q \in \mathbb{Z}.
\end{equation}
We can see that $\{ \varphi_q\}_{q \in\mathbb Z}$ is a partition of unity in $\R^n.$  

Denoting the Fourier transform and its inverse by $\mathcal{F}$ and $\mathcal{F}^{-1},$ respectively, we introduce $h:=\mathcal F^{-1}\varphi.$ For $u \in \mathcal{S}',$ the homogeneous Littlewood-Paley projections are defined as 
\begin{equation}\notag
\dot \Delta_qu:=\mathcal F^{-1}(\varphi(2^{-q}\xi)\mathcal Fu)=2^{nq}\displaystyle\int_{\R^n} h(2^qy)u(x-y)dy, \ q \in \mathbb{Z}.
\end{equation}

In view of the above definitions, we note that the following identity holds in the sense of distributions - 
$$ u= \sum_{q \in \mathbb Z} \dot \Delta_q u.$$
With each $\dot \Delta_q u$ supported in some annular domain in the Fourier space, Littlewood-Paley projections  provide us with a way to decompose a function into pieces with localized frequencies.

For $s \in \R$ and $1 \leq p,q \leq \infty,$ we define the homogeneous Besov space $\dot B^s_{p,q}$ as 
$$\dot B^s_{p,q}(\R^n) = \big\{f \in \mathcal{S}'(\R^n) : \|f\|_{\dot B^{s}_{p,q}(\R^n)} < \infty \big \},$$ with the norm given by
\begin{equation}\notag
\|f\|_{\dot B^{s}_{p,q}(\R^n)}=
\begin{cases}
\displaystyle\big(\sum_{j \in \mathbb Z} (2^{sj}\|\dot \Delta_j f \|_{L^p(\R^n)})^q \big)^{\frac{1}{q}}, \ \text{ if }1 \leq q < \infty,\\
\displaystyle \sup_{j \in \mathbb Z} (2^{sj}\|\dot \Delta_j f \|_{L^p(\R^n)}), \ \text{ if } q= \infty.
\end{cases}
\end{equation}
In this paper, we are primarily interested in the $L^\infty, \ell^\infty$-based Besov spaces $\dot B^{s}_{\infty, \infty}.$

\subsection{Besov spaces and the heat kernel}

It turns out that negative order Besov spaces can also be characterized via the action of the heat kernel. In particular, we have the following lemma, for whose proof we refer readers to \cite{Lr}.
\begin{Lemma} Let $f \in {\dot B^{s}_{\infty, \infty}}$ for some $s<0.$ The following norm equivalence holds.
\begin{equation}
\|f\|_{\dot B^{s}_{\infty, \infty}}=\sup_{t>0} t^{-\frac{s}{2\a}}\|e^{-t(-\Delta)^{\a}}f\|_{L^\infty}, \text{ where } \a >0.
\end{equation}
\end{Lemma}

More generally, the following lemma concerning the action of the heat semigroup in Besov spaces holds true and shall be extensively used in this paper.
\begin{Lemma}\label{bsv}
i) For $\a>0,$ the following inequalities hold.
\begin{gather*}
\|e^{-t(-\Delta)^{\a}}f\|_{L^\infty} \leq  C \|f\|_{L^\infty},\\
\|\nabla e^{-t(-\Delta)^{\a}}f\|_{L^\infty} \leq  C t^{-\frac{1}{2\a}}\|f\|_{L^\infty},\\
\|\nabla \mathbb{P}e^{-t(-\Delta)^{\a}}f\|_{L^\infty} \leq  C t^{-\frac{1}{2\a}}\|f\|_{L^\infty}.
\end{gather*}
ii) For $\a>0$ and $s_0 \leq s_1,$ the following inequalities hold.
\begin{gather*}
\|e^{-t(-\Delta)^{\a}}f\|_{\dot B^{s_1}_{\infty, \infty}} \leq  C t^{-\frac{1}{2\a}(s_1-s_0)} \|f\|_{\dot B^{s_0}_{\infty, \infty}},\\
\|\nabla^k e^{-t(-\Delta)^{\a}}f\|_{\dot B^{s_1}_{\infty, \infty}} \leq  C t^{-\frac{1}{2\a}(s_1-s_0+k)} \|f\|_{\dot B^{s_0}_{\infty, \infty}}.
\end{gather*}
\end{Lemma}
Proofs of Lemma \ref{bsv} can be found in \cite{KOT, MYZ}.

\subsection{Mild solutions}

A mild solution to system (\ref{fhmhd}) is the fix point of the map
\begin{equation}\label{msl}
S(u,b):=\begin{pmatrix}S_1(u,b) \\ S_2(u,b)\end{pmatrix},
\end{equation}
where $S_1(u,b)$ and $S_2(u,b)$ are given by the following Duhamel's formulae -
\begin{equation}\label{msu} 
\begin{split}
S_1(u,b):= u(t,x)=& e^{-\nu t (-\Delta)^{\a_1}}u_0(x)- \int_0^t e^{-\nu (t-s)(-\Delta)^{\a_1}}\mathbb{P}\nabla \cdot (u \otimes u)(s) \mathrm{d}s\\
& + \int_0^t e^{-\nu (t-s)(-\Delta)^{\a_1}}\mathbb{P}\nabla \cdot (b \otimes b)(s) \mathrm{d}s,\\
\end{split}
\end{equation}
\begin{equation}\label{msb}
\begin{split}
S_2 (u,b):= b(t,x) =& e^{-\mu t (-\Delta)^{\a_2}}b_0(x)- \int_0^t e^{-\mu (t-s)(-\Delta)^{\a_2}}\mathbb{P}\nabla \cdot (u \otimes b)(s) \mathrm{d}s\\
& + \int_0^t e^{-\mu (t-s)(-\Delta)^{\a_2}}\mathbb{P}\nabla \cdot (b \otimes u)(s) \mathrm{d}s\\
& - \eta \int_0^t e^{-\mu (t-s)(-\Delta)^{\a_2}}\nabla \times( \nabla \cdot (b \otimes b))(s) \mathrm{d}s.\\
\end{split}
\end{equation}
In (\ref{msb}), we have applied the vector identity $\nabla \times (\nabla \cdot (b \otimes b))=\nabla \times ((\nabla \times b)\times b)$ to the Hall term. To further simplify notations, we view the integrals in expressions (\ref{msu}) and (\ref{msb}) as bilinear forms.
\begin{Definition}[Bilinear forms] Let $f,g \in \mathcal{S}^{'}.$ The bilinear forms $\mathcal{B}_{\a_1}(\cdot,\cdot),$ $\mathcal{B}_{\a_2}(\cdot, \cdot)$ and $\mathfrak{B}_{\a_2}(\cdot, \cdot)$ are defined as follows.
\begin{equation}\notag
\begin{split}
\mathcal{B}_{\a_1}(f,g)=& \int_0^t e^{-\nu (t-s)(-\Delta)^{\a_1}}\mathbb{P}\nabla \cdot (f \otimes g)(s) \mathrm{d}s;\\
\mathcal{B}_{\a_2}(f,g)=& \int_0^t e^{-\mu (t-s)(-\Delta)^{\a_2}}\mathbb{P}\nabla \cdot (f \otimes g)(s) \mathrm{d}s;\\
\mathfrak{B}_{\a_2}(f,g)= & \eta \int_0^t e^{-\mu (t-s)(-\Delta)^{\a_2}}\nabla \times( \nabla \cdot (b \otimes b))(s) \mathrm{d}s.
\end{split}
\end{equation}
\end{Definition}
In view of the above, we can write the formulae (\ref{msl}), (\ref{msu}) and (\ref{msb}) as
\begin{equation}
\begin{split}
S_1(u,b)= & \tilde u_0(x) - \mathcal{B}_{\a_1} (u,u) + \mathcal{B}_{\a_1} (b,b),\\
S_2(u,b)= & \tilde b_0(x) - \mathcal{B}_{\a_2} (u,b) + \mathcal{B}_{\a_2} (b,u) - \mathfrak{B}_{\a_2} (b,b).\\
\end{split}
\end{equation}

\subsection{The contraction principle}

Given the mild solution formulation (\ref{msl}), a traditional approach is to find a fixed point by iterating the map $(u,b) \mapsto S(u,b).$ In order to do so, it is essential to find a space $\mathcal{E}$ such that the bilinear forms $\mathcal{B}_{\a}(\cdot, \cdot)$ and $\mathfrak{B}_{\a}(\cdot, \cdot)$ are bounded from $\mathcal{E} \times \mathcal{E}$ to $\mathcal{E}.$ In this paper, we shall use the following lemma, proven in \cite{Lr} and \cite{M} as a simple consequence of Banach fixed point theorem.
\begin{Lemma}\label{pic}
Let $\mathcal{E}$ be a Banach space. Given a bilinear form $\mathbb{B}: \mathcal{E} \times \mathcal{E} \to \mathcal{E}$ such that $\|\mathbb{B}(u,v)\|_\mathcal{E} \leq C_0 \|u\|_\mathcal{E} \|v\|_\mathcal{E}, \forall u, v \in \mathcal{E},$ for some constant $C_0 >0,$ we have the following assertions for the equation 
\begin{equation}\label{eq1}
u= y +\mathbb{B}(u,u).
\end{equation}

i). Suppose that $y \in B_\varepsilon(0):= \{f\in \mathcal{E}: \|f\|_{\mathcal{E}}< \varepsilon \}$ for some $ \varepsilon \in \big(0, \frac{1}{4C_0}\big),$ then the equation (\ref{eq1}) has a solution $u \in B_{2\varepsilon}(0):= \{f \in \mathcal{E}: \|f\|_{\mathcal{E}}<2\varepsilon \},$ which is, in fact, the unique solution in the ball  $\overline{B_{2\varepsilon}(0)}.$

ii). On top of i), suppose that $\bar y \in B_\varepsilon(0), \bar u \in B_{2\varepsilon}(0)$ and $\bar u=\bar y+\mathbb{B}(\bar u, \bar u),$ then the following continuous dependence is true.
\begin{equation}\label{inq}
\|u-\bar u\|_\mathcal{E} \leq \frac{1}{1-4\varepsilon C_0}\|y -\bar y\|_\mathcal{E}.
\end{equation}
\end{Lemma}
It can be seen from inequality (\ref{inq}) that to ensure local well-posedness, it suffices that $C_0=CT^a$ for some $a>0,$ while global well-posedness would require $C_0$ to be bounded above by a time-independent constant. 

\bigskip

\section{Proofs of Theorems}

This section is devoted to the proofs of Theorems \ref{lwp} and \ref{gwp}. We work within a framework based on the concepts of the ``admissible path space" and ``adapted value space", as formulated in \cite{Lr}. The idea is to first identify an ``admissible path space" $\mathcal{E}_T$ in which we may apply the contraction principle, then characterize the ``adapted value space" $E_T$ associated with $\mathcal{E}_T.$ In our case, we consider the space $$
E_T =\{f: f \in \mathcal{S}', \ e^{-t(-\Delta)^{\a_i}}f \in \mathcal{E}_T, \ 0< t <T\}, \ i=1 \text{ or } 2.
$$

To start, we define the Banach spaces $X_T$ and $Y_T$ and the admissible path space $\mathcal{E}_T:=X_T \times Y_T.$
\begin{equation}\label{exy}
X_T=\Big\{f: \R^{+} \to L^\infty(\R^3): \nabla \cdot f=0 \text{ and }\sup_{0< t <T} t^{\frac{2\a_1-\gamma}{2\a_1}}\|f(t)\|_{L^\infty(\R^3)} <\infty \Big \}
\end{equation}
\begin{align}\label{eyx}
Y_T=\Big\{f: \R^{+} \to L^\infty(\R^3): \nabla \cdot f=0 \text{ and }\sup_{ 0< t <T} t^{\frac{2\a_2-\beta}{2\a_2}}\|f(t)\|_{L^\infty(\R^3)} <\infty \Big \}
\end{align}

By formulae (\ref{msu}) and (\ref{msb}) along with the characterization of homogeneous Besov spaces in terms of the heat flow (\ref{bsv}), we have the following inequalities - 
\begin{equation}\notag
\begin{split}
\|u\|_{X_T} \leq & \sup_{t>0} t^{\frac{2\a_1-\gamma}{2\a_1}}\| \tilde u_0 \|_\infty + \|\mathcal{B}_{\a_1}(u,u)\|_{X_T} + \|\mathcal{B}_{\a_1}(b,b)\|_{X_T}\\
\leq & C_\nu \|u_0 \|_{ \dot B^{-(2\a_1-\gamma)}_{\infty, \infty} } + \|\mathcal{B}_{\a_1}(u,u)\|_{X_T} + \|\mathcal{B}_{\a_1}(b,b)\|_{X_T},
\end{split}
\end{equation}
\begin{equation}\notag
\begin{split}
\|b\|_{Y} \leq & \sup_{t>0} t^{\frac{2\a_2-\beta}{2\a_2}}\| \tilde b_0 \|_\infty + \|\mathcal{B}_{\a_2}(u,b)\|_{Y_T}+ \|\mathcal{B}_{\a_2}(b,u)\|_{Y_T}+\|\mathfrak{B}_{\a_2}(b,b)\|_{Y_T}\\
\leq & C_\mu \|b_0 \|_{ \dot B^{-(2\a_2-\beta)}_{\infty, \infty} }+ \|\mathcal{B}_{\a_2}(u,b)\|_{Y_T}+ \|\mathcal{B}_{\a_2}(b,u)\|_{Y_T}+\|\mathfrak{B}_{\a_2}(b,b)\|_{Y_T}.
\end{split}
\end{equation}
Clearly, ${\dot B^{-(2\a_1-\gamma)}_{\infty, \infty}} \times {\dot B^{-(2\a_2-\b)}_{\infty, \infty}(\R^3)}$ is an adapted value space corresponding to the admissible path space $\mathcal{E}_T$ given by Definitions \ref{exy} and \ref{eyx}. 

We proceed to prove the following proposition.

\begin{Proposition}
Suppose that the parameters $\a_1, \a_2, \b$ and $\gamma$ satisfy 
\begin{equation}\label{cnt}
\begin{cases}
\gamma \geq \max\{1, \frac{\a_1}{\a_2}\}, \\
\b \geq \max\{ 2, \frac{(\gamma +1)\a_2}{2\a_1}\},\\
\frac{\gamma}{2}< \a_1 < \gamma,\\
\frac{\b}{2} <\a_2<\b.
\end{cases} 
\end{equation}
If $(u,b) \in \mathcal{E}_T$ for some $0<T<\infty,$ then $\|S(u,b)-(\tilde u_0, \tilde b_0)\| \in \mathcal{E}_T.$ In particular, 
\begin{equation}\label{bbd}
\|S(u,b)-(\tilde u_0, \tilde b_0)\|_{\mathcal{E}_T} \leq C T^a \|(u,b)\|_{\mathcal{E}_T}^2
\end{equation}
for some $a>0$ and $C=C(\nu, \mu, \eta)>0.$
\end{Proposition}

\pf First, we remark that the constraints on the parameters indeed yield a non-empty set, since the combination $\gamma=1, \beta=2, \a_1=1-\delta$ and $\a_2=2-2\delta$ with $\frac{1}{4}< \delta <  \frac{1}{2}$ clearly satisfies (\ref{cnt}). 

To prove (\ref{bbd}), it suffices to show that the bilinear forms are bounded from $\mathcal{E}_T\times \mathcal{E}_T$ to $\mathcal{E}_T,$ with bounds dependent on $\nu, \mu, \eta$ and $T.$ To this end, we invoke the property of the Beta function. More specifically, for $\a>1$ and $0 < \theta < \a,$ we have
\begin{equation}\label{bet}
\int^t_0 (t-\tau)^{-\frac{1}{\a}}\tau^{-\frac{\theta}{\a}}\mathrm{d}\tau = t^{1-\frac{1}{\a}-\frac{\theta}{\a}}B\Big(1-\frac{\theta}{\a}, 1-\frac{1}{\a}\Big) \leq C t^{1-\frac{1}{\a}-\frac{\theta}{\a}}.
\end{equation}

Let $\gamma \geq 1$ and $\frac{\gamma}{2} < \a_1 < \gamma.$ Via integration by parts, H\"older's inequality, identity (\ref{bet}) and Definition \ref{exy}, we have the following inequalities.  
\begin{equation}\notag
\begin{split}
\|\mathcal{B}_{\a_1}(u,u)\|_{X_T} \leq & C_\nu \sup_{0<t<T} t^{\frac{2\a_1-\gamma}{2\a_1}}\int^t_0(t-s)^{-\frac{1}{2\a_1}}\|u(s)\|_\infty\|u(s)\|_\infty \mathrm{d}s\\
\leq & C_\nu \|u\|_{X_T}^2 \sup_{0<t<T} t^{\frac{2\a_1-\gamma}{2\a_1}}\int^t_0(t-s)^{-\frac{1}{2\a_1}}s^{-2+\frac{\gamma}{\a_1}} \mathrm{d}s\\
\leq & C_\nu T^{\frac{\gamma-1}{2\a_1}} \|u\|_{X_T}^2.
\end{split}
\end{equation}

Similarly, the following estimates are true provided that $\gamma \geq 1$, $\frac{\gamma}{2} < \a_1 < \gamma$, $\frac{\beta}{2} < \a_2 < \beta$ and $\b \geq \frac{(\gamma+1)\a_2}{2\a_1}.$
\begin{equation}\notag
\begin{split}
\|\mathcal{B}_{\a_1}(b,b)\|_{X_T} \leq & C_\nu \sup_{0<t<T} t^{\frac{2\a_1-\gamma}{2\a_1}}\int^t_0(t-s)^{-\frac{1}{2\a_1}}\|b(s)\|_\infty\|b(s)\|_\infty \mathrm{d}s\\
\leq & C_\nu \|b\|_{X_T}^2\sup_{0<t<T} t^{\frac{2\a_1-\gamma}{2\a_1}}\int^t_0(t-s)^{-\frac{1}{2\a_1}}s^{-2+\frac{\beta}{\a_2}} \mathrm{d}s\\
\leq & C_\nu T^{\frac{\beta}{\a_2}-\frac{\gamma+1}{2\a_1}} \|b\|_{X_T}^2.
\end{split}
\end{equation}

To bound the term $\|\mathcal{B}_{\a_2}(b,u)\|_Y,$ we further require that $\a_2 > \frac{1}{2}$ and $\gamma \geq \frac{\a_1}{\a_2}.$ 
\begin{equation}\notag
\begin{split}
\|\mathcal{B}_{\a_2}(b,u)\|_{Y_T} \leq & C_\mu \sup_{0<t<T} t^{\frac{2\a_2-\beta}{2\a_2}}\int^t_0(t-s)^{-\frac{1}{2\a_2}}\|u(s)\|_\infty\|b(s)\|_\infty \mathrm{d}s\\
\leq & C_\mu \|u\|_X \|b\|_Y \sup_{0<t<T} t^{\frac{2\a_2-\beta}{2\a_2}}\int^t_0(t-s)^{-\frac{1}{2\a_2}}s^{-2+\frac{\gamma}{2\a_1}+\frac{\beta}{2\a_2}}\mathrm{d}s\\
\leq & C_\mu T^{\frac{\gamma}{2\a_1}-\frac{1}{2\a_2}} \|u\|_{X_T} \|b\|_{Y_T}.
\end{split}
\end{equation}
We note that the term $\|\mathcal{B}_{\a_2}(u,b)\|_Y$ can be estimated in an identical manner. 

Finally, we integrate by parts twice to estimate the Hall term. We end up with the condition $\a_2>1$ along with all the constraints from previous estimates.
\begin{equation}\notag
\begin{split}
\|\mathfrak{B}_{\a_2}(b,b)\|_{Y_T} \leq & C_{\mu, \eta} \sup_{0<t<T} t^{\frac{2\a_2-\beta}{2\a_2}}\int^t_0(t-s)^{-\frac{1}{\a_2}}\|b(s)\|_\infty\|b(s)\|_\infty \mathrm{d}s\\
\leq & C_{\mu, \eta} \|b\|_{Y_T}^2 \sup_{0<t<T} t^{\frac{2\a_2-\beta}{2\a_2}}\int^t_0(t-s)^{-\frac{1}{\a_2}}s^{-2+\frac{\beta}{\a_2}}\mathrm{d}s \\
\leq & C_{\mu, \eta} T^{\frac{\beta-2}{2\a_2}}\|b\|_{Y_T}^2.
\end{split}
\end{equation}

\cbdu

\textbf{Proof of Theorem \ref{lwp}}: By inequality (\ref{bbd}), Lemma \ref{bsv} and Lemma \ref{pic}, there exists a solution $(u,b) \in \mathcal{E}_T$ provided that the initial data $(u_0, b_0)$ and the time $T$ satisfy 
\begin{equation}\notag
4 C T^a \Big(C_{\nu} \|u_0 \|_{ \dot B^{-(2\a_1-\gamma)}_{\infty, \infty} }+  C_{\mu, \eta}\|b_0 \|_{ \dot B^{-(2\a_2-\beta)}_{\infty, \infty} }\Big) < 1.
\end{equation}

It remains to be shown that $(u,b) \in L^\infty\big(0,T; {\dot B^{-(2\a_1-\gamma)}_{\infty, \infty}} \times {\dot B^{-(2\a_2-\b)}_{\infty, \infty}(\R^3)}\big).$ By (\ref{msu}) and Lemma \ref{bsv}, it holds that 
\begin{equation}\notag
\begin{split}
\|S_1 u(t)\|_{\dot B^{-(2\a_1-\gamma)}_{\infty, \infty}} = & \sup_{0<\tau <T}\tau^{\frac{2\a_1-\gamma}{2\a_1}} \|e^{-\nu \tau (-\Delta)^{\a_1}}S_1 u(t) \|_{L^\infty}\\
\lesssim & \sup_{0<\tau <T} \tau^{\frac{2\a_1-\gamma}{2\a_1}} \|e^{-\nu(\tau+t)(-\Delta)^{\a_1}}u_0 \|_{L^\infty}\\
& +\sup_{0<\tau <T} \tau^{\frac{2\a_1-\gamma}{2\a_1}}\|u\|_{X_T}^2\int^{\tau+t}_0 (\tau+t-s)^{-\frac{1}{2\a_1}}s^{-2+\frac{\gamma}{\a_1}}\mathrm{d}s \\
& +\sup_{0<\tau <T} \tau^{\frac{2\a_1-\gamma}{2\a_1}}\|b\|_{Y_T}^2\int^{\tau+t}_0 (\tau+t-s)^{-\frac{1}{2\a_1}}s^{-2+\frac{\b}{\a_2}}\mathrm{d}s.
\end{split}
\end{equation} 
Estimating with the help of (\ref{bet}), we have
\begin{equation}\notag
\begin{split}
\|S_1 u(t)\|_{\dot B^{-(2\a_1-\gamma)}_{\infty, \infty}} \lesssim & \sup_{0<\tau <T} \tau^{\frac{2\a_1-\gamma}{2\a_1}} \Big(\|e^{-\nu\tau(-\Delta)^{\a_1}}u_0 \|_{L^\infty} +(\tau+t)^{-1+\frac{2\gamma-1}{2\a_1}}\|u\|_{X_T}^2\\
&+ (\tau+t)^{-1-\frac{1}{2\a_1}+\frac{2\b}{2\a_2}}\|b\|_{Y_T}^2 \Big)\\
\lesssim & \|u_0\|_{\dot B^{-(2\a_1-\gamma)}_{\infty, \infty} }+ T^a\|(u,b)\|_{\mathcal{E}_T}^2.
\end{split}
\end{equation}

In a similar fashion, the following inequalities follows from (\ref{msb}) and Lemma \ref{bsv}. 
\begin{equation}\notag
\begin{split}
\|S_2 b(t)\|_{\dot B^{-(2\a_2-\b)}_{\infty, \infty}} = & \sup_{0<\tau <T}\tau^{\frac{2\a_2-\b}{2\a_2}} \|e^{-\mu \tau (-\Delta)^{\a_2}}S_2 b(t) \|_{L^\infty}\\
\lesssim & \sup_{0<\tau <T} \tau^{\frac{2\a_2-\b}{2\a_2}} \bigg( \|e^{-\mu(\tau+t)(-\Delta)^{\a_2}}b_0 \|_{L^\infty}\\
& + 2 \|u\|_{X_T}\|b\|_{Y_T}\int^{\tau+t}_0 (\tau+t-s)^{-\frac{1}{2\a_2}}s^{-2+\frac{\gamma}{2\a_1}+\frac{\b}{2\a_2}}\mathrm{d}s \\
& +\|b\|_{Y_T}^2\int^{\tau+t}_0 (\tau+t-s)^{-\frac{1}{\a_2}}s^{ -2+\frac{\b}{\a_2}}\mathrm{d}s \bigg).
\end{split}
\end{equation} 
The integrals can be evaluated thanks to (\ref{bet}), which yields the bound on $S_2 b.$
\begin{equation}\notag
\begin{split}
\|S_2 b(t)\|_{\dot B^{-(2\a_2-\b)}_{\infty, \infty}} \lesssim & \sup_{0<\tau <T} \tau^{\frac{2\a_2-\b}{2\a_2}} \Big(\|e^{-\nu\tau(-\Delta)^{\a_1}}b_0 \|_{L^\infty}\\
& +(\tau+t)^{-1+\frac{\gamma}{2\a_1}+\frac{\b-1}{2\a_2}}\|u\|_{X_T} \|b\|_{Y_T}+ (\tau+t)^{-1+\frac{\b-1}{\a_2}}\|b\|_{Y_T}^2 \Big)\\
\lesssim & \|b_0\|_{\dot B^{-(2\a_2-\b)}_{\infty, \infty}}+ T^a\|(u,b)\|_{\mathcal{E}_T}^2.
\end{split}
\end{equation}
The inequalities above imply that $$(u,b) \in L^\infty\big(0,T; {\dot B^{-(2\a_1-\gamma)}_{\infty, \infty}} \times {\dot B^{-(2\a_2-\b)}_{\infty, \infty}(\R^3)}\big).$$

\cbdu

However,  well-posedness result for the standard Hall-MHD system, i.e., the case $\a_1=\a_2=1,$ is unattainable as the above method breaks down in this case.

We now turn to the hyper-resistive EMHD equations, written as
\begin{equation}\label{emhde}
\begin{cases}
b_t + \eta \nabla \times ((\nabla \times b) \times b)= -\mu (-\Delta)^{\a_2}b,\\
\nabla\cdot b=0,\\
b(0,x)=b_0, \ t \in \R^+, x \in \R^3,
\end{cases}
\end{equation}
where $1< \a_2 <2.$

The above system is the small-scale limit of the Hall-MHD system, corresponding to the scenario in which the ions are practically static, simply forming a neutralizing background for the moving electrons. It is named electron MHD as the system is solely determined by the electrons. In astrophysics, system (\ref{emhde}) makes frequent appearances in the study of the magnetosphere and the solar wind, whose dynamics can be puzzling due to high frequency magnetic fluctuations. Readers may consult \cite{G1, G3, MG} for relevant physics backgrounds.

Unlike the complete system (\ref{fhmhd}),  system (\ref{emhde}) possesses the property of scaling invariance. More specifically, if $b(t,x)$ solves system (\ref{emhde}) with initial data $b_0,$ then $b_\lambda(t,x)= \lambda^{2\a_2-2}b(\lambda^{2\a_2}t, \lambda x)$ is a solution subject to the initial data $\lambda^{2\a_2-2}b_0(\lambda x).$ One can see that the space $L^\infty\big(0, \infty; \dot B^{-(2\a_2-2)}_{\infty, \infty}(\R^3)\big)$ is the largest critical space according to the scaling property. 

We proceed to prove Theorem \ref{gwp} by finding a ball $B \subset Y_T$ where the solution map $S_2$ is a contraction mapping. We have the following two propositions.

\begin{Proposition} Let $\a_2 \in (1,2)$ and $\beta =2.$ For $0< T \leq \infty,$ the map $S_2$ satisfies
\begin{equation}\label{bd2}
\|S_2 b - \tilde b_0\|_{Y_T} \leq C \|b\|_{Y_T}^2.
\end{equation}
Therefore, there exists some $\varepsilon_1 >0,$ such that $S_2$ is a self-mapping on the ball $$B_{\varepsilon_1}\big(\tilde b_0\big)=:\{ f\in Y_T: \|f-\tilde b_0\|_{Y_T}  < \varepsilon_1 \},$$ provided that $ \|b_0\|_{\dot B^{-(2\a_2-2)}_{\infty, \infty}(\R^3)} <\varepsilon_1.$
\end{Proposition}

\pf The inequality (\ref{bd2}) follows from the following estimate.
\begin{equation}\notag
\begin{split}
\|\mathfrak{B}_{\a_2}(b,b)\|_{Y_T} \leq & \sup_{t>0} t^{\frac{2\a_2-2}{2\a}}\int^t_0(t-s)^{-\frac{1}{\a_2}}\|b(s)\|_\infty\|b(s)\|_\infty \mathrm{d}s\\
\leq & \|b\|_{Y_T}^2 \sup_{t>0} t^{\frac{2\a_2-2}{2\a_2}}\int^t_0(t-s)^{-\frac{1}{\a_2}}s^{-2+\frac{2}{\a_2}}\mathrm{d}s \\
\leq & C_{\mu, \eta} \|b\|_{Y_T}^2.
\end{split}
\end{equation}

Since it is assumed that $b \in B_{\varepsilon_1}\big(\tilde b_0\big)$ and $\|b_0\|_{\dot B^{-(2\a_2-2)}_{\infty, \infty}(\R^3)} <\varepsilon_1,$ it follows from inequality (\ref{bd2}) and lemma  (\ref{bsv}) that
\begin{equation}\notag
\|S_2 b - \tilde b_0\|_{Y_T} \leq C \|b\|_{Y_T}^2
\leq C \big( \|b-\tilde b_0\|_{Y_T}^2 + \|\tilde b_0\|_{Y_T}^2 \big)
\leq C \varepsilon_1^2.
\end{equation}

\cbdu

\begin{Proposition}\label{ctr}
Let $1< a_2 <2$ and $\beta =2.$ For any $T \in (0, \infty],$ there exists some $\varepsilon_2 \in (0, \varepsilon_1)$ such that if $\|b_0\|_{\dot B^{-(2\a_2-2)}_{\infty, \infty}(\R^3)} <\varepsilon_2,$ then the solution map $S_2$ is a contraction mapping on the ball $$B_{\varepsilon_2}\big(\tilde b_0\big)=:\{ f\in Y_T: \|f-\tilde b_0\|_{Y_T}  < \varepsilon_2 \}.$$
\end{Proposition}

\pf Let $b, \bar b \in B_{\varepsilon_2}\big( \tilde b_0 \big).$ Clearly, the following inequalities hold.
\begin{equation}\notag
\begin{split}
\|S_2 b -S_2 \bar b \|_{Y_T}= & \| \mathfrak{B}_{\a_2}(b,b) - \mathfrak{B}_{\a_2}(\bar b, \bar b) \|_{Y_T}\\
\leq & \| \mathfrak{B}_{\a_2}(b,b) - \mathfrak{B}_{\a_2}(b, \bar b)\|_{Y_T} +  \| \mathfrak{B}_{\a_2}(b, \bar b)  - \mathfrak{B}_{\a_2}(\bar b, \bar b) \|_{Y_T}\\
\leq & C_{\mu, \eta} \max \{\|b\|_{Y_T}, \|\bar b \|_{Y_T}\} \|b-\bar b \|_{Y_T}\\
\leq & C_{\mu, \eta} \varepsilon_2 \|b-\bar b \|_{Y_T}.
\end{split}
\end{equation}
We can ensure that $S_2$ is a contraction mapping by choosing $\varepsilon_2< 1/{2C_{\mu, \eta}}.$

\cbdu

\textbf{Proof of Theorem \ref{gwp}}. As a result of Proposition \ref{ctr}, we know that for some $\varepsilon_2>0,$ $S_2$ has a fixed point, which is a mild solution to system (\ref{emhde}), in $$B_{\varepsilon_2}\big(\tilde b_0\big)=:\{ f\in Y_T: \|f-\tilde b_0\|_{Y_T}  < \varepsilon_2, \ T=+\infty\},$$ provided that $ \|b_0\|_{\dot B^{-(2\a_2-2)}_{\infty, \infty}(\R^3)} <\varepsilon_2.$

To see that the solution $b$ is in $L^\infty(0, \infty; \dot B^{-(2\a_2-2)}_{\infty, \infty}(\R^3)),$ we just calculate
\begin{equation}\notag
\begin{split}
\|S_2 b(t)\|_{\dot B^{-(2\a_2-2)}_{\infty, \infty}} \lesssim & \sup_{\tau >0} \tau^{\frac{2\a_2-2}{2\a_2}} \bigg( \|e^{-\mu(\tau+t)(-\Delta)^{\a_2}}b_0 \|_{L^\infty}\\
& + \|b\|_{Y_T}^2\int^{\tau+t}_0 (\tau+t-s)^{-\frac{1}{\a_2}}s^{ -2+\frac{2}{\a_2}}\mathrm{d}s \bigg)\\
\lesssim & \|b_0\|_{ \dot B^{-(2\a_2-2)}_{\infty, \infty} }+\|b\|_{Y_T}^2.
\end{split}
\end{equation}

\cbdu

Unfortunately, the above pathway to small data global well-posedness fails just when $\a_2=1,$ leaving the question of the standard EMHD equations' solvability in its largest critical space $\dot B^0_{\infty, \infty}(\R^3)$ unanswered. At this moment, we are inclined to believe that in this setting, the system is ill-posed instead.

\bigskip

{\textbf{Acknowledgement.}}
The authors would like to thank Prof. Isabelle Gallagher and Dr. Trevor Leslie for helpful discussions.


\begin{thebibliography}{XX}


\bibitem{ADFL}
M. Acheritogaray, P. Degond, A. Frouvelle and J. Liu.
\newblock {\em Kinetic formulation and global existence for the Hall-Magneto-hydrodynamics system}.
\newblock Kinet. Relat. Models Vol. 4(4), 901-918, 2011.

\bibitem{BF}
M. J. Benvenutti and L. C. F. Ferreira.
\newblock {\em  Existence and stability of global large strong solutions for the Hall-MHD system.}
\newblock Differ. Integral Equ. Vol. 29(9–10), 977–1000, 2016.

\bibitem{BCD}
H. Bahouri, J. Chemin,  and R. Danchin.
\newblock {\em Fourier Analysis and Nonlinear Partial Differential Equations}.
\newblock Grundlehren der mathematischen Wissenschaften, 343. Springer, Heidelberg, 2011.

\bibitem{C}
L. M. B. C. Campos.
\newblock {\em On hydromagnetic waves in atmospheres with application to the Sun.}
\newblock Theor. Comput. Fluid Dyn. 10 (1-4), 37-70, 1998.


\bibitem{CDL}
D. Chae, P. Degond and J. Liu.
\newblock {\em  Well-posedness for Hall-magnetohydrodynamics}.
\newblock Ann. Inst. H. Poincar\'e Anal. Non Lin\'eaire Vol. 31, No. 3, 555-565, 2014.

\bibitem{CL}
D. Chae and J. Lee.
\newblock {\em On the blow-up criterion and small data global existence for the Hall-magnetohydrodynamics.}
\newblock J. Diff. Eq. Vol. 256(11), 3835-3858, 2014.

\bibitem{CS}
D. Chae and M. E. Schonbek.
\newblock {\em On the temporal decay for the Hall-magnetohydrodynamic equations.}
\newblock J. Diff. Eq. Vol. 255(11), 3971-3982, 2013.

\bibitem{CWW}
D. Chae, R. Wan and J. Wu.
\newblock {\em  Local well-posedness for the Hall-MHD equations with fractional magnetic diffusion.}
\newblock J. Math. Fluid Mech. Vol. 17(4), 627-638, 2015.

\bibitem{CW}
D. Chae and S. Weng.
\newblock {\em  Singularity formation for the incompressible Hall-MHD equations without resistivity.}
\newblock Ann. Inst. H. Poincar\'e Anal. Non Lin\'eaire Vol. 33, No. 4, 1009-1022, 2016.

\bibitem{D1}
M. Dai.
\newblock {\em Regularity criterion for the 3D Hall-magneto-hydrodynamics.}
\newblock J. Diff. Eq. Vol. 261(1), 573-591, 2016.

\bibitem{D2}
M. Dai.
\newblock {\em  Local well-posedness of the Hall-MHD system in $H^s(\R^n)$ with $s>n/2$.}
\newblock arXiv: 1709.02347.

\bibitem{D3}
M. Dai.
\newblock {\em Local well-posedness for the Hall-MHD system in optimal Sobolev Spaces.}
\newblock arXiv: 1803.09556.

\bibitem{D4}
M. Dai. 
\newblock {\em  Non-uniqueness of Leray-Hopf weak solutions of the 3D Hall-MHD system}.
\newblock arXiv: 1812.11311.

\bibitem{DL}
M. Dai and H. Liu.
\newblock {\em Long time behavior of solutions to the 3D Hall-magneto-hydrodynamics system with one diffusion}.
\newblock J. Diff. Eq., Vol. 266, 7658-7677, 2019.

\bibitem{FFNZ}
J. Fan, Y. Fukumoto, G.Nakamura and Y. Zhou.
\newblock {\em Regularity criteria for the incompressible Hall-MHD system  .}
\newblock Z. Angew. Math. Mech. 95(11),  1156-1160, 2015. 

\bibitem{FLN} 
J. Fan, F. Li and G. Nakamura.
\newblock {\em Regularity criteria for the incompressible Hall-magnetohydrodynamic equations .}
\newblock Nonlinear Anal. 109: 173-179, 2014.

\bibitem{FSZ}
J. Fan, B. Samet and Y. Zhou.
\newblock {\em A regularity criterion for a generalized Hall-MHD system.}
\newblock Comput. Math. Appl. 74(10), 2438-2443, 2017.

\bibitem{G1}
S. Galtier.
\newblock{\em Introduction to Modern Magnetohydrodynamics.}  
\newblock Cambridge University Press, Cambridge, UK, 2016.

\bibitem{G2}
S. Galtier.
\newblock {\em Wave turbulence in incompressible Hall magnetohydrodynamics .}
\newblock Journal of Plasma Physics 72 (5), 721-769, 2006.

\bibitem{GB}
S. Galtier and E Buchlin.
\newblock {\em Multiscale Hall-magnetohydrodynamic turbulence in the solar wind .}
\newblock The Astrophysical Journal 656 (1), 560, 2007.

\bibitem{G3}
S. Galtier.
\newblock {\em Exact scaling laws for 3D electron MHD turbulence.}
\newblock Journal of Geophysical Research: Space Physics 113 (A1), 2008.

\bibitem{G}
L. Grafakos.
\newblock {\em Modern Fourier Analysis .}
\newblock Graduate Texts in Mathematics, Vol. 250, 2nd edition, Springer, New York, 2009.

\bibitem{GMS}
W. Gu, C. Ma and J. Sun.
\newblock {\em A regularity criterion for the generalized Hall-MHD system.}
\newblock Bound. Value Probl. 2016:188, 2016.

\bibitem{HAHZ}
F. He, B. Ahmad, T. Hayat and Y. Zhou.
\newblock {\em  On regularity criteria for the 3D Hall-MHD equations in terms of the velocity.}
\newblock Nonlinear Anal. RWA 32, 35-51, 2016.

\bibitem{JO}
I. Jeong and S. Oh.
\newblock{\em On the Cauchy problem for the Hall and electron magnetohydrodynamic equations without resistivity I: illposedness near degenerate stationary solutions.}
\newblock arXiv: 1902.02025.

\bibitem{JZ}
Z. Jiang and M. Zhu.
\newblock {\em  Regularity criteria for the 3D generalized MHD and Hall-MHD systems.}
\newblock Bull. Malays. Math. Sci. Soc. 41 (1), 105-122, 2018.

\bibitem{KOT}
H. Kozono, T. Ogawa and Y. Taniuchi.
\newblock {\em  Navier-Stokes equations in the Besov space near $L^\infty$ and BMO.}
\newblock Kyushu Journal of Mathematics 57:303-324, 2003.

\bibitem{KL}
M. Kwak and B. Lkhagvasuren.
\newblock {\em Global wellposedness for Hall-MHD equations }.
\newblock Nonlinear Anal. 174: 104-117, 2018.

\bibitem{Lr}
P. G. Lemari\'e-Rieusset
\newblock {\em Recent Development in the Navier-Stokes Problem.}
\newblock Chapman \& Hall/CRC Press: Boca Raton, 2002. 

\bibitem{L}
M. J. Lighthill.
\newblock {\em Studies on magneto-hydrodynamic waves and other anisotropic wave motions.}
\newblock Phil. Trans. R. Soc. A
252 (1014), 397-430, 1960.

\bibitem{M}
Y. Meyer.
\newblock {\em Wavelets, paraproducts and Navier-Stokes equations .}
\newblock  Current Development in Mathematics, 1996, p105-212, International Press, Cambridge, MA, 1999.

\bibitem{MG}
R. Meyrand and S. Galtier.
\newblock {\em  Anomalous  Spectrum in Electron Magnetohydrodynamic Turbulence.}
\newblock Physical review letters 111 (26), 264501, 2013.

\bibitem{MYZ}
C. Miao, B. Yuan and B. Zhang.
\newblock {\em  Well-posedness of the Cauchy problem for the fractional power dissipative equations.}
\newblock Nonlinear Analysis: Theory, Methods \& Applications 68: 461-484, 2008.

\bibitem{PMZ}
N. Pan, C. Ma and M. Zhu. 
\newblock {\em Global regularity for the 3D generalized Hall-MHD system.}
\newblock Appl. Math. Lett. 61: 62-66, 2016.

\bibitem{PZ}
N. Pan and M. Zhu.
\newblock {\em A new regularity criterion for the 3D generalized Hall-MHD system with $\b \in (\frac12,1]$.}
\newblock J. Math. Anal. Appl, Vol. 445(1), 604-611, 2017.

\bibitem{PM}
J. M. Polygiannakis and X. Moussas.
\newblock {\em  A review of magneto-vorticity induction in Hall-MHD plasmas.}
\newblock  Plasma Phys. Control. Fusion 43 (2), 195, 2001.

\bibitem{W}
R. Wan. 
\newblock {\em Global regularity for generalized Hall Magneto-Hydrodynamics systems .}
\newblock  Electron. J. Differ. Equ. 179, 1-18, 2015.

\bibitem{WZ1}
R. Wan and Y. Zhou. 
\newblock {\em  On global existence, energy decay and blow-up criteria for the Hall-MHD system.}
\newblock J. Diff.
Eq. 259 (11), 5982-6008, 2015.

\bibitem{WZ2}
R. Wan and Y. Zhou.
\newblock {\em  Low Regularity Well-Posedness for the 3D Generalized Hall-MHD System.}
\newblock Acta Appl. Math. 147, 95-111, 2017.

\bibitem{WL}
Y. Wang and H. Li.
\newblock {\em Beale-Kato-Madja type criteria of smooth solutions to 3D Hall-MHD flows.}
\newblock Appl. Math. Comput. 286, 41-48, 2016.

\bibitem{WYT}
X. Wu, Y. Yu and Y. Tang
\newblock {\em Global existence and asymptotic behavior for the 3D generalized Hall-MHD system.}
\newblock Nonlinear Anal. 151:41-50, 2017.

\bibitem{Y1}
Z. Ye.
\newblock {\em Regularity criteria and small data global existence to the generalized viscous Hall-magnetohydrodynamics.}
\newblock Comput. Math. Appl. 70(8), 2137-2154, 2015.

\bibitem{Y2}
Z. Ye
\newblock {\em  Global well-posedness and decay results to 3D generalized viscous magnetohydrodynamic equations.}
\newblock Ann. Mat. Pura Appl. (4) 195 (4), 1111-1121, 2016. 

\bibitem{Y3}
Z. Ye.
\newblock {\em Regularity criterion for the 3D Hall-magnetohydrodynamic equations involving the vorticity.}
\newblock Nonlinear Anal. 144, 182-193, 2016. 

\bibitem{Y4}
Z. Ye.
\newblock {\em A logarithmically improved regularity criterion for the 3D Hall-MHD equations in Besov spaces with negative 
indices.}
\newblock Appl. Anal. 96 (16), 2669-2683, 2017. 

\bibitem{YZ0}
Z. Ye and Z. Zhang.
\newblock {\em A remark on regularity criterion for the 3D Hall-MHD equations based on the vorticity.}
\newblock Appl. Math. Comput. 301: 70–77, 2017.

\bibitem{YZ}
X. Yu and Z. Zhai.
\newblock {\em Well-posedness for fractional Navier-Stokes equations in the largest critical spaces $\dot B^{-(2\b-1)}_{\infty, \infty}(\R^n)$.}
\newblock Math. Meth. Appl. Sci. 35(6), 676-683, 2012.

\bibitem{Z}
Z. Zhang.
\newblock {\em A remark on the blow-up criterion for the 3D Hall-MHD system in Besov spaces.}
\newblock J. Math. Anal. Appl. 441 (2), 692-701, 2016.

\end{thebibliography}
\end{document}